\documentclass[12pt]{article}
\usepackage{fullpage,graphicx,psfrag,amsmath,amsfonts,amssymb}
\usepackage[small,bf]{caption}
\usepackage{graphicx}
\usepackage{url,algorithm,amsmath,amsthm}
\usepackage{algorithmicx,algpseudocode}
\usepackage{subcaption}
\newcommand{\BEAS}{\begin{eqnarray*}}
\newcommand{\EEAS}{\end{eqnarray*}}
\newcommand{\BEA}{\begin{eqnarray}}
\newcommand{\EEA}{\end{eqnarray}}
\newcommand{\BEQ}{\begin{equation}}
\newcommand{\EEQ}{\end{equation}}
\newcommand{\BIT}{\begin{itemize}}
\newcommand{\EIT}{\end{itemize}}
\newcommand{\BNUM}{\begin{enumerate}}
\newcommand{\ENUM}{\end{enumerate}}

\newcommand{\BA}{\begin{array}}
\newcommand{\EA}{\end{array}}

\newcommand{\eg}{{\it e.g.}}
\newcommand{\ie}{{\it i.e.}}

\newcommand{\ones}{\mathbf 1}

\newcommand{\reals}{{\mbox{\bf R}}}

\newcommand{\symm}{{\mbox{\bf S}}}  

\newcommand{\Tr}{\mathop{\bf Tr}}
\newcommand{\diag}{\mathop{\bf diag}}

\newcommand{\argmin}{\mathop{\rm argmin}}

\newtheorem{theorem}{Theorem}
\newtheorem{lemma}{Lemma}
\newtheorem{assumption}{Assumption}

{\begin{quote}}{\end{quote}}

\makeatletter
\long\def\@makecaption#1#2{
   \vskip 9pt 
   \begin{small}
   \setbox\@tempboxa\hbox{{\bf #1:} #2}
   \ifdim \wd\@tempboxa > 5.5in
        \begin{center}
        \begin{minipage}[t]{5.5in}
        \addtolength{\baselineskip}{-0.95pt}
        {\bf #1:} #2 \par
        \addtolength{\baselineskip}{0.95pt}
        \end{minipage}
        \end{center}
   \else 
	\hbox to\hsize{\hfil\box\@tempboxa\hfil}  
   \fi
   \end{small}\par
}
\makeatother

\newcounter{oursection}

\newcounter{lecture}

\newcommand{\Psucc}{p^\mathrm{succ}}

\title{Linear Programming Heuristics for the Graph Isomorphism Problem}
\author{Reza Takapoui \and Stephen Boyd}
\begin{document}
\maketitle
\date{\vspace{-5ex}}

\begin{abstract}
An isomorphism between two graphs is a bijection between 
their vertices that preserves the edges. 
We consider the problem of determining whether two finite 
undirected weighted graphs are isomorphic, and 
finding an isomorphism relating them if the answer is positive. 
In this paper we introduce effective probabilistic linear programming (LP) 
heuristics to solve 
the graph isomorphism problem. We motivate our heuristics by showing
guarantees under some conditions, and present numerical experiments that 
show effectiveness of these heuristics in the general case.
\end{abstract}

\section{Graph isomorphism problem} 
\subsection{Problem statement}

Consider two weighted undirected graphs, each with $n$ vertices labeled $1, \ldots, n$, 
described by their adjacency matrices $A,\tilde A\in\reals^{n\times n}$, 
where $A_{ij}$ is the weight on the edge in the first graph between vertices $i$ and $j$,
and zero if there is no edge between vertices $i$ and $j$ (and similarly for $\tilde A$).
Since the graphs are undirected, the adjacency matrices are symmetric.
The two graphs are isomorphic if there is a permutation of the vertices of the first graph 
that makes the first graph the same as the second. This occurs if and only if there is a 
permutation matrix $P \in \reals^{n \times n}$ (\ie, a matrix with exactly one entry in each row 
and column that is one, with the other zero) that satisfies $PAP^T=\tilde A$. 
We will say that the permutation matrix $P$ transforms $A$ to $\tilde A$ if 
$PAP^T=\tilde A$, or equivalently, $PA=\tilde AP$.

The graph isomorphism problem (GIP) is to determine whether such a permutation 
matrix exists, and to find one if so.
GIP can be formulated as the (feasibility) optimization problem
 \begin{equation}
    \begin{array}{ll}
    \mbox{find}   &P \\
    \mbox{subject to} & PA=\tilde AP\\
     & P\ones= \ones,\quad P^T\ones= \ones\\
    & P_{ij}\in\{0,1\},\quad i,j=1,\ldots,n,
    \end{array}
    \label{problem}
 \end{equation}
with variable $P\in \reals^{n \times n}$, where $\ones$ is the vector
with all entries one. The data for this problem are the adjacency 
matrices of the two graphs, $A$ and $\tilde A$.
The constraints on the last two lines enforce that $P$ is a permutation matrix. 
This problem has $n^2$
scalar variables (\ie, $P_{ij}$), and all constraints except the last one (which requires
the variables to take on Boolean values $0$ or $1$) are linear in $P$.
The problem~(\ref{problem}) is an integer (or Boolean) linear program (LP).

If $P$ transforms $A$ to $\tilde A$, the two matrices are similar, and therefore have the same 
spectrum, \ie, the same eigenvalues including multiplicity. This observation gives a very 
simple spectral condition for isomorphism: The spectra of $A$ and $\tilde A$ must be the 
same. Since this condition is readily checked, we will henceforth assume that this is the case.

\begin{assumption}\label{a-spectra}
Adjacency matrices $A$ and $\tilde A$ have the same eigenvalues including multiplicity.
\end{assumption}
Assumption \ref{a-spectra} does not imply that the graphs are isomorphic. 
But if assumption~\ref{a-spectra} does not hold, the graphs are surely not isomorphic.

GIP arises in several applications
including chem-informatics, mathematical chemistry, electronic design automation,
and network theory. For example it can be used in 
determining if two descriptions of a molecule are the same \cite{irniger2005graph}, or whether 
the physical layout of an electronic circuit correctly reflects the 
given circuit schematic diagram \cite{colbourn1981testing}.

\subsection{Computational complexity}
GIP is one of the few problems in NP that has so far 
resisted all attempts to be classified as 
NP-complete, or within P. In $1979$, Garey and Johnson \cite{garey1979computers}
listed $12$ such problems and the GIP is one of only two on that list 
whose complexity remains unresolved.
For many restricted graph classes, polynomial time
algorithms are known. This is, for example, 
the case for trees \cite{kelly1957congruence},
planar graphs \cite{hopcroft1974linear},
interval graphs \cite{lueker1979linear},
partial $k$-trees \cite{bodlaender1990polynomial},
graphs of bounded degree \cite{luks1980isomorphism},
or graphs with bounded eigenvalue multiplicity \cite{babai1982isomorphism}.
Recently, Babai announced a quasipolynomial time algorithm for all graphs,
\ie, one with running time $2^{O\left((\log n)^c\right)}$ for some fixed
$c>0$ \cite{babai2015graph}.

This paper does not advance the discussion of computational complexity
of GIP;
instead we describe effective heuristics for it.  However,
many of the ideas we will encounter (\eg, 
compact graphs, graph spectrum) are closely connected to those
that arise in research on its computational complexity.

\subsection{A simple relaxation}\label{s:rlx}
One relaxation of problem~(\ref{problem}) can be obtained by 
replacing constraints $P_{ij}\in\{0,1\}$ with interval
relaxation $0\leq P_{ij}\leq1$. 
The relaxed problem can be written as
\begin{equation}
    \begin{array}{ll}
    \mbox{find}   &P \\
    \mbox{subject to} & PA=\tilde AP\\
     & P\ones= \ones,\quad P^T\ones= \ones\\
    & P \geq 0,
    \end{array}
    \label{simple-rlx}
 \end{equation}
 where the last constraint is entry-wise. 
 The constraint $P\leq 1$ is removed because it is implied by $P\ones=P^T\ones=\ones$ and 
 $P\geq 0$. This problem has $n^2$ scalar variables, 
$n^2+2n$ scalar linear equality constraints (some of which are
redundant), and $n^2$ scalar linear inequality constraints. 
While~(\ref{problem}) is a Boolean LP, the 
relaxation~(\ref{simple-rlx}) is an LP, and 
readily solved in polynomial time using different 
methods such as interior-point methods
or ellipsoid method \cite{boyd2004convex}.

If the relaxed problem~(\ref{simple-rlx}) is infeasible, 
then $A$ and $\tilde A$ are obviously not isomorphic.
Therefore, from this point of
the paper on, we assume that this problem is feasible.
\begin{assumption}\label{a-feasible}
We assume the relaxed problem~(\ref{simple-rlx}) is feasible.
\end{assumption}
Assumption \ref{a-feasible} does not imply that the graphs are isomorphic. 
But if assumption~\ref{a-feasible} does not hold, the graphs are surely not isomorphic.
The necessary and sufficient conditions for feasibility of this 
LP has been
studied in \cite{scheinerman2011fractional} and is 
closely related to common equitable partitions of two graphs. 

The set of feasible points for problem~(\ref{simple-rlx}) 
is a polytope in $\reals^{n\times n}$. 
Permutation matrices transforming $A$ to $\tilde A$ are
extreme points or vertices of this polytope. (However, not every 
extreme point of this polytope is necessarily a
permutation matrix.) 
Given assumption~\ref{a-feasible}, GIP
comes down to finding a permutation matrix among extreme points of
the feasible set of problem~(\ref{simple-rlx}).
Our algorithms are all based on variations and extensions of
this observation. 

Surprisingly, under some conditions on the graphs 
(stated in \S\ref{theory}),
this relaxation is tight; that is,
the set of feasible point for problems~(\ref{problem}) and (\ref{simple-rlx})
are both singletons. Hence,
it is sufficient to solve problem~(\ref{simple-rlx}) to find a 
permutation
matrix transforming $A$ to $\tilde A$,
or certify that $A$ and $\tilde A$ are not isomorphic
if problem~(\ref{simple-rlx}) is infeasible.

\subsection{Outline}
We will describe the basic version of the algorithm in 
\S\ref{algorithm}, then describe sparsity
constraints in \S\ref{improve} to tighten the relaxation and increase the probability 
of success for the heuristic.
An ADMM-based algorithm to solve the problem is described in \S\ref{solve}, and finally
in \S\ref{examples} we study some examples and numerical experiments. 

\section{Algorithm}\label{algorithm}
\subsection{Basic randomized subroutine}\label{basic}
As we discussed in \S\ref{s:rlx}, the set of feasible points for problem~(\ref{simple-rlx})
form a polytope and permutation matrices transforming $A$ to $\tilde A$ are
extreme points or vertices of this polytope.
In order to encourage 
finding an extreme point of this polytope, we minimize 
a random linear objective over the set of feasible points.
Hence our basic randomized subroutine is to generate a random matrix
$W\in\reals^{n\times n}$ with i.i.d.\ Gaussian entries and 
solve the LP
 \begin{equation}
    \begin{array}{ll}
    \mbox{minimize}   &\Tr(W^TP) \\
    \mbox{subject to} & PA=\tilde AP\\
    & P\ones= \ones, \quad P^T\ones= \ones\\
    &P \geq 0.
    \end{array}
    \label{heuristic}
 \end{equation}

If the solution to problem~(\ref{heuristic}) found 
happens to be a permutation matrix
we conclude that we solved this instance of GIP, since
$A$ and $\tilde A$ are isomorphic and the solution $P^\star$ is a permutation 
matrix that transforms $A$ to $\tilde A$.
On the other hand, if the solution found is not a permutation matrix, 
we cannot say anything about $A$ and $\tilde A$. 

Notice that with probability $1$ the solution of 
problem~(\ref{heuristic}) is unique.
Let $\Psucc$ denote the probability under the distribution of $W$ that the
randomized subroutine finds a permutation matrix when
$A$ and $\tilde A$ are isomorphic. It is easy to see that 
$\Psucc$ is only a function of $A$
(and does not depend on $\tilde A$). We briefly comment here that
$\Psucc>0$.

Since scaling $W$ does not change the optimization problem,
without loss of generality we can assume that $W$ is chosen
uniformly from an $(n^2-1)$-sphere in an $n^2$ dimensional space.
Consider the set of permutation matrices in $\reals^{n\times n}$
that transform $A$ to $\tilde A$. 
Centered around each such permutation matrix, there is 
a cone with nonzero solid angle such that if $-W$ is in that cone,
the solution of problem~(\ref{heuristic}) 
will be a permutation matrix. Probability of success
$\Psucc$ is the probability of $W$ being in one of these cones, 
and hence the probability that this algorithm successfully finds
an acceptable permutation matrix. 

Even though solving just one instance of problem~(\ref{heuristic})
might look unsophisticated, 
it turns out to be an effective heuristic
to solve the GIP in some cases.  It is guaranteed to find a
solution of GIP (with probability one) 
under some conditions which are discussed in next subsection. 

Standard solvers based on the simplex method and interior-point methods
can be used to solve problem~(\ref{heuristic}), but can be inefficient
in practice.  However, we will show
in \S\ref{solve} that the special structure of this problem enables us to use a custom 
method to solve this optimization problem more efficiently.

\subsection{Polishing}\label{polish}
After the basic randomized subroutine converges 
to a doubly stochastic matrix $P^*\in\reals^{n\times n}$, 
we use can the Hungarian algorithm \cite{kuhn1955hungarian} to project $P^*$ on the
set of doubly stochastic matrices, which we denote by $\tilde P^*$. Finding the closest
permutation matrix to $P^*$ (in Frobenius norm) can be done in $O(n^3)$ time.
If $\tilde P^*A=\tilde A\tilde P^*$, then
a permutation mapping $A$ to $\tilde A$ is found. This step, as will see in
\S\ref{examples}, can increase the probability
of success of the subroutine.

\subsection{Repeated randomized algorithm}
As we discussed in \S\ref{basic}, choosing a random instance of $W$ 
and running the basic randomized algorithm will find a permutation that relates
$A$ and $\tilde A$ with probability $\Psucc$. 
We also discussed in \S\ref{polish} that we can use Hungarian algorithm 
to (potentially) increase the probability of success.
If we repeat this basic randomized
algorithm for $N$ independently chosen random 
instances of $W$ the probability of success is 
$1-(1-\Psucc)^N$. For example if the probability of success for the basic randomized
subroutine is $80\%$, after repeating the basic subroutine $4$ times the probability of
success will be $99.9\%$. 

We evidently have $\Psucc>0$, meaning that by 
solving problem~(\ref{heuristic}) there is a positive probability of
finding a permutation matrix that that transforms $A$ to
$\tilde A$.  In particular, by repeatedly solving  
problem~(\ref{heuristic}) with independent random choices of $W$,
we will solve GIP with probability one.
This probability can be extremely small, however, so
the expected number of LPs we must solve to solve GIP can be 
extremely large.
The benefit of this probabilistic method over a deterministic spectral-based 
heuristic is that when a deterministic heuristic fails, there is no hope to 
recover the permutation
that relates $A$ and $\tilde A$, but this probabilistic method can be used repeatedly to 
find the underlying permutation.
This randomized algorithm can potentially have false negatives, \ie, 
when $A$ and $\tilde A$ are
isomorphic, it might fail to find a permutation relating them.

\subsection{Theoretical guarantees}\label{theory}
In this subsection we show that under some 
conditions, solving problem~(\ref{heuristic}) is 
guaranteed to find a permutation solution if $A$ 
and $\tilde A$ are isomorphic, in other words, $\Psucc = 1$.
Specifically
if $A$ has distinct eigenvalues, and
for every eigenvector $v$ of $A$ we have $\ones^Tv\neq 0$, the relaxation in 
problem~(\ref{simple-rlx}) is tight.

In order to show this, assume $QA=\tilde AQ$, where $Q$ is a permutation matrix and
let  $P$ be a doubly stochastic matrix such that $PA=\tilde AP$. Also let
$A=V\Lambda V^T$ be an eigendecompostion of $A$. Defining 
$R = V^TQ^TPV$, we have
\[
R\Lambda = V^TQ^TPV \Lambda = V^TQ^TPAV = V^TQ^T\tilde A PV = V^TAQ^T PV =
\Lambda V^TQ^TPV = \Lambda R.
\]
Hence for every $i,j$ we have $R_{ij} \Lambda_{ii} = \Lambda_{jj} R_{ij}$, or equivalently
$R_{ij}\left( \Lambda_{ii} - \Lambda_{jj} \right) = 0$. Since
$\Lambda_{ii}-\Lambda_{jj}\neq 0$ for $i\neq j$, all off-diagonal entries of $R$ must be
zero and the matrix $R$ must be diagonal. Also
$RV^T\ones = V^TQ^TP\ones = V^TQ^T\ones = V^T\ones$, which implies that
$R_{ii} (V^T\ones)_i = (V^T\ones)_i$. Our second assumption on $V$ enforces
that $(V^T\ones)_i\neq 0$ and
$R_{ii} = 1$ for all $i$ and hence $P = QVRV^T = QVV^T = Q$.

For graphs specified above, relaxing the set of permutation matrices to the set of
doubly stochastic matrices does not extend the set of feasible points 
of problem~(\ref{problem}), which is a singleton.

A graph with adjacency matrix $A$ is said to be a \emph{compact graph} if
the set of feasible points to problem~(\ref{heuristic}) is the convex hull of 
the set of feasible points to problem~(\ref{problem}).
For example, all graphs with two aforementioned properties are compact.
The concept of \emph{compact graphs} was introduced by Tinhofer 
\cite{tinhofer1991note}, who 
proved that trees and cycles (which violate the two 
aforementioned assumptions) and the disjoint union of 
isomorphic copies of a compact graph are compact. 
If $A$ is the adjacency matrix of a compact graph and $A$ and $\tilde A$
are isomorphic, then $\Psucc=1$. 
It is not a hard problem to solve the graph isomorphism 
problem between two trees
or two cycles (or some other compact graphs), but it is interesting 
to see that this problem can be solved (with probability 
$1$) by only solving one LP.

\section{Sparsity constraints}\label{improve}

We discussed one possible relaxation of problem~(\ref{problem}) 
in \S\ref{s:rlx}. However,
we can find tighter relaxations by adding convex constraints 
about $P$ that we know 
must be true. This will create a tighter relaxation of the original 
problem which can potentially increase $\Psucc$ for the randomized
algorithm. Specifically, we consider the following extension 
 \begin{equation}
    \begin{array}{ll}
    \mbox{minimize}   &\Tr(W^TP) \\
    \mbox{subject to} & PA=\tilde AP\\
    & P\ones= \ones, \quad P^T\ones= \ones\\
    &P \geq 0\\
    &P_{ij} = 0 \quad (i,j)\in \mathcal K,
    \end{array}
    \label{improved}
 \end{equation}
where $\mathcal K$ is a set of pairs of indices. 
The difference between this extension 
and problem~(\ref{heuristic}) is the last line of the constraints,
which requires some entries in $P$ to be zero.
If $\mathcal K=\emptyset$, this problem reduces to
problem~(\ref{heuristic}). In general, problem~(\ref{improved})
can be considered as a problem with $n^2-\|\mathcal K\|$
scalar variables. Our efficient method for solving the LP relaxation
(described in \S\ref{solve}) can handle these constraints
efficiently.

There are different ways to find a proper set $\mathcal K$. The simplest
and least useful choice of $\mathcal K$ is the empty set. The maximal
$\mathcal K$, denoted by $\mathcal K^{\max}$ is the set of pairs of indices
$i,j$ such that $P_{ij}=0$ for every permutation matrix $P$ transforming 
$A$ to $\tilde A$. Any valid set $\mathcal K$
will satisfy $\emptyset \subseteq \mathcal K \subseteq \mathcal K^{\max}$.

We will show below how to find pairs of 
indices $i,j$ such that for any permutation matrix $P$ satisfying
$PA=\tilde A P$, we must have $P_{ij}=0$.

\begin{lemma}
Let $P$ be a permutation matrix. If $Pa=b$ for two given vectors 
$a,b\in\reals^n$ and $P_{ij}=1$, then we must have $a_j=b_i$. This implies 
that if $a_j\neq b_i$ then $P_{ij}=0$. 
\end{lemma}
This simple lemma can be used to construct $\mathcal K$.
If we know that $Pa=b$ for any permutation matrix
that transforms $A$ to $\tilde A$, we can include the pair $(i,j)$ in
$\mathcal K$.
If all entries of vector $a$ are distinct, then $\mathcal K$ contains 
$n^2-n$ entries and permutation $P$
is uniquely determined. 
At the other extreme, if all entries of vector $a$ are equal, this lemma adds no
new information.

The idea here is to use graph invariants to find equalities of the form 
$Pa=b$ and then
construct a set $\mathcal K$ from these equalities.
By graph invariant, we mean 
any function of nodes that is independent
of labeling the nodes. For instance, the number of paths 
with length $k$ starting from a node
is an invariant, where $k$ is a positive integer. 
(We will discuss this in next subsection.) 
In principle, the method can be used with any 
invariant.  But we are interested in invariants that are 
easy to compute, \eg,
degree invariants and spectral invariants.

\subsection{Degree invariants}
A simple linear equality that holds is $P(\diag A) = \diag \tilde A$. This means that 
vertices that are mapped to each other must have the same self-edge
weights. If $A$ and $\tilde A$ have no self-edges
(\ie, diagonal entries of $A$ and $\tilde A$ are zero), this equality will add no information.

A more interesting equality that holds is $P(A^k\ones) = \tilde A^k \ones$ 
for every positive integer $k$.
In other words, the number of paths of length $k$ starting from a vertex
must be the same for two vertices that map to each other. 

For example, when the graphs are unweighted 
(\ie, the edges have weight one)
this equality
with $k=1$ means that the degrees of mapped vertices must be the same. 
In other words, all nodes of degree $i$ must be mapped to each other. Hence,
if the number of nodes with degree $i$ is denoted by $n_i$, the original problem
has $(n_1+n_2+\cdots)^2$ variables, but we know only $n_1^2+n_2^2+\cdots$
of them can potentially be nonzero. In the extreme 
case that all degrees are equal (or equivalently the graph is regular), this
method does not eliminate any variables from $P$.

This equality is valid for every positive integer $k$. However, according to 
the Cayley-Hamilton theorem, $A^k$ can be written as a linear combination of
$A,A^2,\ldots,A^n$ for every positive $k$. Hence for a given pair of indices $i,j$ if
$(A^k)_i=(A^k)_j$ for $k=1,\ldots,n$, we will have $(A^k)_i=(A^k)_j$
for every positive $k$. Therefore it is only enough to consider this equality
for $k=1,\ldots,n$.

In summary, we have two sets of equalities:
\BIT
\item $P(\diag A) = \diag \tilde A$.
\item $P(A^k\ones) = \tilde A^k \ones$, for $k=1,\ldots,n$. 
\EIT
We denote the set of pairs of indices constructed this way by 
$\mathcal K^\mathrm{degree}$. Clearly, we have
$\mathcal K^\mathrm{degree}\subseteq\mathcal K^{\max}$.

\subsection{Spectral invariants}
As mentioned earlier, we assume that $A$ and 
$\tilde A$ share the same spectrum; otherwise the graphs are
evidently non-isomorphic.
Let $\lambda$ be an eigenvalue of $A$ with multiplicity $k$ and 
columns of $V\in\reals^{n\times k}$ be an orthonormal basis for 
eigenspace associated with $\lambda$, hence 
we have $AV = \lambda V$ and $V^TV = I_k$, where $I_k$ denotes the
identity matrix in $\reals^{k\times k}$. Similarly, assume that
columns of $\tilde V\in\reals^{n\times k}$ are an orthonormal basis of eigenspace
associated with $\lambda$, hence 
we have $\tilde A\tilde V = \lambda \tilde V$ and $\tilde V^T\tilde V= I_k$.
We have
\[
\tilde A(PV) = (\tilde AP)V = (PA) V = P(AV) = P (\lambda V) = \lambda (PV).
\]
In other words columns of $PV$ are eigenvectors of $\tilde A$ associated with $\lambda$.
Therefore we have $PV = \tilde V Q$ where $Q\in\reals^{k\times k}$. We observe that
\[
Q^TQ = Q^T\tilde V^T\tilde VQ = V^TP^TPV = V^TV = I_k.
\]
Hence $Q$ is an orthogonal in $\reals^{k\times k}$. Therefore $PV=\tilde VQ$ implies that
$P(VV^T)P^T = \tilde V \tilde V^T$, and we have the following equalities:
\BIT
\item $P(\diag VV^T) = \diag \tilde V\tilde V^T.$
\item $P(VV^T\ones) = \tilde V\tilde V^T\ones.$
\EIT
Notice that $P\left(\left(VV^T\right)^k\ones\right)= (\tilde V\tilde V^T)^k\ones$ adds
no more information for $k>1$, since $(VV^T)^k=VV^T$ and 
$(\tilde V\tilde V^T)^k=\tilde V\tilde V^T$.

These equalities hold for any eigenvalue $\lambda$. 
We denote the set of pairs of indices constructed this way by 
$\mathcal K^\mathrm{spectral}$. Clearly, we have
$\mathcal K^\mathrm{spectral}\subseteq\mathcal K^{\max}$.

\subsection{Constructing $\mathcal K$}\label{construct}
In summary, here is how we construct 
$\mathcal K = \mathcal K^\mathrm{degree} \cup \mathcal K^\mathrm{spectral}$. 
We start with $\mathcal K=\emptyset$, and we
sequentially consider the degree invariant equalities for $k=1,\ldots,n$ 
and spectral invariant equalities for
every eigenvalue $\lambda$. For each equality, 
we add the new pairs of disallowed $i,j$ to the set 
$\mathcal K$. 

When $\mathcal K$ is constructed, the number of variables that could possibly be nonzero 
is $n^2-|\mathcal K|$. In \S\ref{examples} we will report the number of elements in $\mathcal K$ for 
our experiments, and we see how the sparsity constraints can 
increase the probability of success for our heuristic, for each LP solved.

Simple examples show that $\mathcal K^\mathrm{degree}$
and $\mathcal K^\mathrm{spectral}$ are not necessarily subsets of each other.
Therefore in general it is a good idea to include both of these constraints
in the basic subroutine. Also, our examples show that \emph{pruning}
can be a helpful technique for making $\mathcal K$ larger. In pruning,
we disallow pair $i,j$ if no neighbor of $i$ can be mapped to a neighbor
of $j$. 

One reasonable conjecture could be that choosing the maximal $\mathcal K$ 
would result in solving the problem with probability $1$. In other words
if $\mathcal K=\mathcal K^{\max}$ then $\Psucc=1$. This conjecture is wrong,
however. It can be shown that for the Petersen graph $\mathcal K^{\max}=\emptyset$.
Also $(1/3)A$ is a doubly stochastic matrix that commutes with $A$.
If this conjecture were true, then $A$ could be 
written as the convex combination of automorphisms of the Petersen graph. 
Each such automorphism $\Pi$ has the property that $\Pi(v)$ is connected 
to $v$ for all vertices $v$. But the only automorphism of the Petersen 
graph that has this property is the identity. Hence the Petersen
graph disproves this conjecture.

\section{Algorithm summary}
A summary of the algorithm is presented below.
\begin{algorithm}[H]
\caption{Randomized LP heuristic} 
\begin{algorithmic}[1]\label{alg}
\State Check whether assumptions~\ref{a-spectra},~\ref{a-feasible} hold.
If not, declare $A$ and $\tilde A$ as non-isomorphic.
\State Construct $\mathcal K$ as described in \S\ref{construct}
\For{random instance $W_1,\ldots,W_N$}
\State Solve problem~(\ref{improved}) to find a solution $P^*$
\State Use Hungarian algorithm to find the closest permutation matrix
$\tilde P^*$ to $P^*$
\If {$\tilde P^*A=\tilde A \tilde P^*$}
\State Declare $A$ and $\tilde A$ as isomorphic and return $\tilde P^*$
\EndIf
\EndFor
\end{algorithmic}
\label{alg_summary}
\end{algorithm}
\bigskip
Here $N$ denotes the number of times problem~(\ref{improved})
is solved and can be chosen beforehand based on an approximation
of $\Psucc$. As mentioned before, the probability of the success of 
this algorithm is $1-(1-\Psucc)^N$.

\section{Solution method} \label{solve}
\subsection{ADMM algorithm}
Any LP solver such as ones base on the simplex method or 
interior-point methods can be used to solve 
problem~(\ref{improved}).
With $n^2$ variables and $n^2$ equality constraints, 
using a general purpose interior-point
solver would involve $O(n^6)$ flops. However, 
the special structure of the problem enables us
to use a custom method which is more efficient.
We use the alternating direction method of multipliers (ADMM) 
in consensus form \cite{boyd2011distributed} 
for solving problem~(\ref{improved}) as follows. We write the problem as
\begin{equation}
    \begin{array}{ll}
    \mbox{minimize}   & f_1(P_1) + f_2(P_2) + f_3(P_3)  \\
    \mbox{subject to} & Z = P_i, \quad i=1,2,3\\
    &ZA = \tilde AZ,
    \end{array}
    \label{consensus}
 \end{equation}
where
\BEAS
f_1(P_1) &=& (1/2) \Tr(W^TP_1) + \mathbf I(P_1\ones = \ones),\\
f_2(P_2) &=& (1/2) \Tr(W^TP_2) + \mathbf I(P_2^T\ones = \ones),\\
f_3(P_3) &=&  \mathbf I(P_3 \geq 0) + \mathbf I
\left(\left(P_3\right)_{i,j} = 0,\quad (i,j)\in\mathcal K\right).
\EEAS
Here $\mathbf I$ is the indicator function, and takes the value $0$
if its argument
is true, the value $\infty$ if its argument is false.
Assuming that $A=V\Lambda V^T$ and $\tilde A=\tilde V\Lambda \tilde V^T$ are
eigenvalue decompositions of $A$ and $\tilde A$,
the ADMM updates will be
\BEAS
P_1^{k+1} &=& \argmin_{P_1\ones = \ones} \left(\Tr\left(\left(W/2+Y_1^k\right)^TW\right)+(1/2)\|P_1-Z^k\|_F^2\right)\\
P_2^{k+1} &=& \argmin_{P_2^T\ones = \ones} \left(\Tr\left(\left(W/2+Y_2^k\right)^TW\right)+(1/2)\|P_2-Z^k\|_F^2\right)\\
P_3^{k+1} &=& \Pi_{\{P|f_3(P)=0\}}(Z^k - Y_3^k)\\
Z^{k+1} &=& \Pi_{\{Z|ZA=\tilde AZ\}}\left(\frac{P_1^{k+1}+P_2^{k+1}+P_3^{k+1}}3\right)\\
Y_1^{k+1} &=& Y_1^k + P_1^{k+1} - Z^{k+1}\\
Y_2^{k+1} &= &Y_2^k + P_2^{k+1} - Z^{k+1}\\
Y_3^{k+1} &= &Y_3^k + P_3^{k+1} - Z^{k+1},
\EEAS
where $\Pi_\mathcal C$ denotes the projection operator onto $\mathcal C$.
Without loss of generality, we chose $\rho=1$ in ADMM, because scaling $\rho$ is 
equivalent to scaling $W$. After generating a random direction $W$, we scale $W$
such that the absolute value of elements of $W$ averages to $1$. 
This is only for the ADMM algorithm to converge faster, and it does not
affect the solution found by the ADMM. (Remember that the solution
to problem~(\ref{consensus}) is unique with probability $1$.)

In order to find the projection onto $\{Z|ZA=\tilde AZ\}$ we notice that $ZA=\tilde AZ$ is
equivalent to $ZV\Lambda V^T =  \tilde V\Lambda \tilde V^T Z$, which is 
equivalent to $\tilde V^T ZV\Lambda =  \Lambda \tilde V^T Z V^T$, or equivalently
$\tilde V^T ZV$ commutes with $\Lambda$. Hence, $ZA=\tilde AZ$ if and only if
$(\tilde V^T ZV)_{ij} = 0$ for any $i,j$ that $\Lambda_{ii} \neq \Lambda_{jj}$.
After simplifying the steps above, the ADMM update steps will reduce to
\BEAS
P_1^{k+1} &=& Z^k - W/2 - Y_1^k - (1/n) \left( \left(Z^k - W/2 - Y_1^k\right) \ones - \ones \right)\ones^T\\
P_2^{k+1} &=& Z^k - W/2 - Y_2^k - (1/n) \ones \left(\ones^T \left(Z^k - W/2 - Y_2^k\right) - \ones^T \right)\\
P_3^{k+1} &=&\max \left( Z^k - Y_3^k, 0 \right)\circ S \\
Z^{k+1} &=& \tilde V \left( \left( \tilde V^T\frac{P_1^{k+1}+P_2^{k+1}+P_3^{k+1}}3V\right)\circ R\right) V^T\\
Y_1^{k+1} &=& Y_1^k + P_1^{k+1} - Z^{k+1}\\
Y_2^{k+1} &= &Y_2^k + P_2^{k+1} - Z^{k+1}\\
Y_3^{k+1} &= &Y_3^k + P_3^{k+1} - Z^{k+1},
\EEAS 
where $\circ$ denotes entry wise product, and $R_{ij}$ is $1$ if $\lambda_i=\lambda_j$ and $0$
otherwise, and $S_{ij} = 1$ if $(i,j)\in K$ and $0$ otherwise.

The first three update rules are the first proximal operator in ADMM algorithm and 
are, in fact, minimizing a quadratic function over an affine subspace. Hence, they are
linear operators in $P$, and are $O(n^2)$ flops.
The fourth equality is the projection of $(P_1^{k+1}+P_2^{k+1}+P_3^{k+1})/3$ onto 
$\{Z|ZA=\tilde AZ\}$, and is of $O(n^3)$ complexity. 
Finally, the last three equalities are dual updates,
and can be done in $O(n^2)$ flops. Therefore each iteration of the ADMM algorithm is done
in $O(n^3)$ flops.

The practical complexity of this heuristic is $O(n^3)$, which is basically the square root of 
the one obtained with interior-point methods. We observe that this algorithm usually converges 
within tens to hundreds of ADMM iterations. 

\section{Examples}\label{examples}
We use several test cases to examine our heuristic. In the first subsection 
we discuss a specific example, the Frucht Graph, in order to give a qualitative description of the 
heuristic. Through this example, we will describe how this heuristic can benefit from
using graph invariant constraints.
Then we use several standard classes of graphs for benchmarking algorithms. We will
report our experiment results for undirected random graphs, two dimensional grids, cubic 
Hype-Hamiltonian graphs, Kronecker Eye Filp graphs, Miayzaky Augmented graphs,
and affine geometries graphs
\cite{de2003large, junttila2011conflict,junttila2007engineering,
lopez2011conauto,lopez2009fast,Specht2013}.
These graphs are described in \S\ref{benchmark}. 
We observe that this heuristic is more effective for some classes of graphs. We also notice that
the time complexity of this heuristic is $O(n^3)$, and it is more effective when graph
invariant constraints are taken into account. The code for our solver and our examples
can be found at \url{https://github.com/cvxgrp/graph_isom}.

\subsection{Frucht graph}
\begin{figure}[h]
\begin{center}

\includegraphics[width=0.25\textwidth]{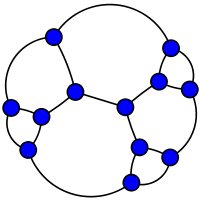}
\end{center}
\caption{\footnotesize The Frucht graph is an asymmetric $3$-regular graph
with $12$ vertices and $18$ edges.}
\label{frucht_img}
\end{figure}

The Frucht graph is an example of a graph with no nontrivial automorphisms,
that violates the assumptions in \S\ref{theory}. As shown in Figure~\ref{frucht_img}, the Frucht 
graph is a $3$-regular graph with $12$ vertices and
$18$ edges. It was first described by Robert Frucht in $1939$. 

We permute nodes of the Frucht graph randomly and 
run the algorithm for $100,000$ times with and without sparsity mask.
Each run uses a random instance of linear objectives ($W$). We observe that
without using sparsity extension, in $25.2\%$ of random instances the algorithm finds the 
permutation matrix. When we used sparsity mask (graph invariant constraints, or $\mathcal K$), 
which contains only $14$ nonzeros,
the algorithm found the correct permutation in every single instance of $100,000$ runs,
so our guess is that $\Psucc=1$ when sparsity constraints are taken into account.

\subsection{Database of graphs for benchmarking algorithms}\label{benchmark}
We use several classes of graphs in \cite{de2003large, junttila2011conflict,junttila2007engineering,
lopez2011conauto,lopez2009fast} to evaluate the effectiveness of this
heuristic on a variety of graphs. In particular, our experiments 
were on the following classes of graphs.
\BIT
\item Undirected random graphs (R1N family) is taken from \cite{de2003large}. In these
undirected Erdos-Renyi graphs, there is an edge between two vertices with probability $0.1$.
\item Cubic Hypo-Hamiltonian graphs clique-connected (CCH family) is proposed
in \cite{lopez2011conauto}, and is built using as basic block two non-isomorphic HypoHamiltonian
graphs with 22 vertices. 
\item Affine geometries graphs (AG2 family) is a part of the benchmark of bliss \cite{Specht2013}.
It contains point-line graphs of affine geometries $AG_2(q)$.
\item Two dimensional grids (G2N family) also comes from \cite{de2003large} benchmark and 
is discussed in \cite{kandel2007applied}. They are introduced 
for simulating applications dealing with regular 
structures as those operating at the lower levels of an image processing task.
\item Kronecker Eye Flip graphs (KEF family) is a part of the benchmark of bliss \cite{Specht2013}.
\item Miayzaky Augmented graphs (MZA family) are also taken from 
the benchmark of bliss \cite{Specht2013}.
\EIT 

For each graph, we ran the algorithm for $50$ random instances 
of $W$ and have summarized the results in Table~\ref{tab}
below. The first column states the graph type, and the second column states the ratio
of possibly nonzero entries in $P$. In other words the number in the second column
is $1-\frac{|\mathcal K|}{n^2}$. The third and fourth column represent the percentage of
runs that a permutation was successfully found with and without sparsity constraints ($\mathcal K$),
respectively.

{\footnotesize

\begin{table}

\begin{center}
 \begin{tabular}{||p{8mm} c p{10mm} p{15mm} p{15mm}||}
 \hline
 {\scriptsize Graph type} & {\scriptsize $n$} & {\scriptsize sparsity ratio} & {\scriptsize success rate without $\mathcal K$} & {\scriptsize success rate with $\mathcal K$}  \\ [0.5ex] 
 \hline\hline
 R1N & 20 & 0.050 & 1.0 & 1.0\\
 \hline
 R1N & 40 & 0.025 & 1.0 & 1.0\\
 \hline
 R1N & 60 & 0.017 & 1.0 & 1.0\\
 \hline
 R1N & 80 & 0.012 & 1.0 & 1.0\\
 \hline
 R1N & 100 & 0.010 & 1.0 & 1.0\\
 \hline
 R1N & 200 & 0.005 & 1.0 & 1.0\\
 \hline
 R1N & 400 & 0.003 & 1.0 & 1.0\\
 \hline
 R1N & 600 & 0.002 & 1.0 & 1.0\\
 \hline
 R1N & 800 & 0.001 & 1.0 & 1.0\\
 \hline
 R1N & 1000 & 0.001 & 1.0 & 1.0\\
 \hline \hline
  CHH & 22 & 0.393 & 0.74 & 0.88\\
 \hline
 CHH & 44 & 0.205 & 0.86 &0.96  \\
 \hline
 CHH & 88 & 0.204 & 0.72 &0.96  \\
 \hline
 CHH & 132 & 0.226 & 0.78 &0.92  \\
 \hline
 CHH & 198 & 0.225 & 0.66  &0.96  \\
 \hline
 CHH & 264 & 0.204 &0.54  &0.90  \\
 \hline
 CHH & 352 & 0.205 & 0.48 &0.86  \\
 \hline
 CHH & 440 & 0.212 &0.32 &0.92  \\
 \hline
 CHH & 550 & 0.212 &0.30 &0.82  \\
 \hline
 CHH & 660 & 0.204 & 0.38 &0.82  \\
 \hline
 CHH & 792 & 0.205 & 0.34 &0.82  \\
 \hline
 CHH & 924 & 0.208 & 0.26 &0.82  \\
 \hline\hline
AG2 & 10 & 0.520 & 0.34 & 1.0  \\
 \hline
AG2 & 21 & 0.510 & 0.20 & 1.0  \\
 \hline
AG2 & 35 & 0.506 & 0.12 & 0.72  \\
 \hline
AG2 & 55 & 0.504 & 0 & 0  \\
 \hline
AG2 & 105 & 0.502 & 0 & 0  \\
 \hline
AG2 & 136 & 0.501 & 0 & 0  \\
 \hline
 \end{tabular}
\begin{tabular}{||p{8mm} c p{10mm} p{15mm} p{15mm}||} 
 \hline
 {\scriptsize Graph type} & {\scriptsize $n$} & {\scriptsize sparsity ratio} & {\scriptsize success rate without $\mathcal K$} & {\scriptsize success rate with $\mathcal K$}  \\ [0.5ex] 
 \hline\hline
  G2N & 16 & 0.375 & 0.54& 0.98 \\
 \hline
 G2N & 36 & 0.186 & 0.52& 0.98 \\
 \hline
 G2N & 64 & 0.110 & 0.46& 1.0 \\
 \hline
 G2N & 81 & 0.078 & 0.44& 0.96 \\
 \hline
 G2N & 100 & 0.072 & 0.42& 1.0 \\
 \hline
 G2N & 196 & 0.038 & 0.58& 0.98 \\
 \hline
 G2N & 400 & 0.019  & 0.52& 0.98\\
 \hline
 G2N & 576 & 0.013 & 0.54& 1.0 \\
 \hline
 G2N & 784 & 0.010 & 0.48& 0.90 \\
\hline\hline
 KEF & 32 & 0.375 & 0.68 & 0.80  \\
 \hline 
 KEF & 50 & 0.270 & 0.82 & 0.86  \\
 \hline 
 KEF & 72 & 0.167 & 0.78 & 0.88  \\
 \hline 
 KEF & 98 & 0.103 & 0.62 & 0.92  \\
 \hline 
 KEF & 128 & 0.062 & 1.0 & 1.0  \\
 \hline 
 KEF & 242 & 0.046 & 1.0 & 1.0  \\
 \hline 
 KEF & 392 & 0.036 & 1.0 & 1.0  \\
 \hline 
 KEF & 578 & 0.029 & 1.0 & 1.0  \\
 \hline 
 KEF & 800 & 0.025 & 1.0 & 1.0  \\
 \hline 
 KEF & 1058 & 0.021 & 1.0 & 1.0  \\
 \hline\hline
 MZA & 40 & 0.159 & 1.0 & 1.0  \\
 \hline
MZA & 80 & 0.085 & 0.62 & 0.76  \\
 \hline
MZA & 120 & 0.057 & 0.22 & 0.28  \\
 \hline
MZA & 160 & 0.044 & 0.24 & 0.32  \\
 \hline
MZA & 200 & 0.035 & 0.24 & 0.32  \\
 \hline
MZA & 240 & 0.029 & 0.08 & 0.16  \\
 \hline
MZA & 280 & 0.025 & 0 & 0.08  \\
 \hline
MZA & 320 & 0.022 & 0 & 0.08  \\
 \hline
MZA & 360 & 0.019 & 0 & 0  \\
 \hline
\end{tabular} 

\end{center}
\caption{The outcome for $50$ instances of our heuristic on different problems.}\label{tab}

\end{table}

}

\bigskip

\begin{figure}[h]
\begin{center}
\includegraphics[width=0.85\textwidth]{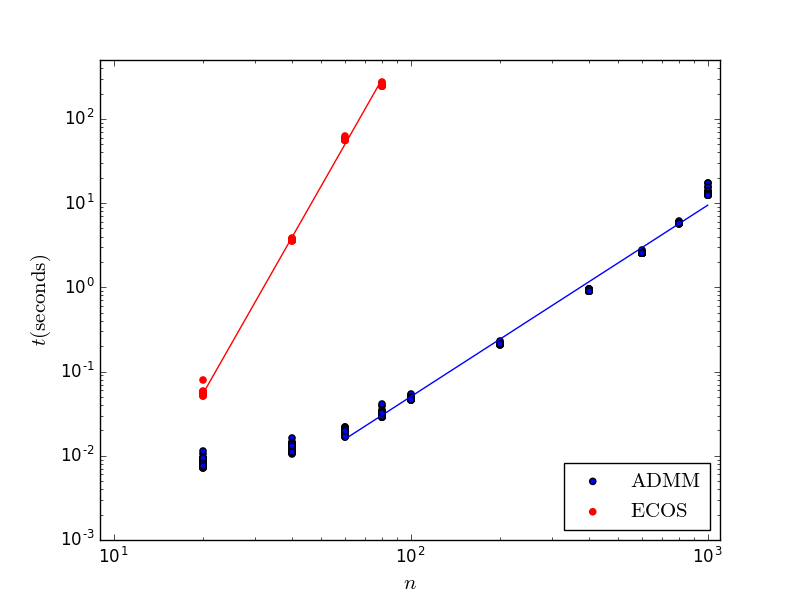}
\end{center}
\caption{\footnotesize Runtime versus graph size in log-scale for 50 random runs.}
\label{time}
\end{figure}

As we discussed in \S\ref{solve}, each iteration of the ADMM algorithm is $O(n^3)$, hence 
we expect that the runtime of the heuristic be a cubic 
function of the number of vertices of the graphs. 
This is because we observe that the number of iterations 
that ADMM requires before it convergence does not change dramatically as the number of 
vertices changes.
Figure~\ref{time} shows the scatter plot of runtime of the algorithm 
(denoted by $t$) versus
the number of vertices (denoted by $n$) in log-scale for $50$ runs over 
random undirected graphs (R1N family). 
This runtime excludes the polishing step, which is done using Hungarian algorithm.
We observe that the average runtime is
linear in the size of the graph when plotted in log-scale. For comparison we also used ECOS 
\cite{domahidi2013ecos}, an interior-point method for second-order cone programming. We
see that our ADMM solver is far more scalable than interior-point methods.

\clearpage
\nocite{*}
\bibliography{graph_isom}

\newcommand{\etalchar}[1]{$^{#1}$}
\begin{thebibliography}{BPC{\etalchar{+}}11}

\bibitem[ABK14]{aflalo2014graph}
Y.~Aflalo, A.~Bronstein, and R.~Kimmel.
\newblock Graph matching: {R}elax or not?
\newblock {\em arXiv preprint arXiv:1401.7623}, 2014.

\bibitem[ABK15]{aflalo2015convex}
Y.~Aflalo, A.~Bronstein, and R.~Kimmel.
\newblock On convex relaxation of graph isomorphism.
\newblock {\em Proceedings of the National Academy of Sciences},
  112(10):2942--2947, 2015.

\bibitem[Bab15]{babai2015graph}
L.~Babai.
\newblock Graph isomorphism in quasipolynomial time.
\newblock {\em arXiv preprint arXiv:1512.03547}, 2015.

\bibitem[BGM82]{babai1982isomorphism}
L.~Babai, D.~Y. Grigorye, and D.~M. Mount.
\newblock Isomorphism of graphs with bounded eigenvalue multiplicity.
\newblock In {\em Proceedings of the fourteenth annual ACM symposium on Theory
  of computing}, pages 310--324. ACM, 1982.

\bibitem[Bod90]{bodlaender1990polynomial}
H.~L. Bodlaender.
\newblock Polynomial algorithms for graph isomorphism and chromatic index on
  partial k-trees.
\newblock {\em Journal of Algorithms}, 11(4):631--643, 1990.

\bibitem[BPC{\etalchar{+}}11]{boyd2011distributed}
S.~Boyd, N.~Parikh, E.~Chu, B.~Peleato, and J.~Eckstein.
\newblock Distributed optimization and statistical learning via the alternating
  direction method of multipliers.
\newblock {\em Foundations and Trends{\textregistered} in Machine Learning},
  3(1):1--122, 2011.

\bibitem[BV04]{boyd2004convex}
S.~Boyd and L.~Vandenberghe.
\newblock {\em Convex Optimization}.
\newblock Cambridge University Press, 2004.

\bibitem[Col81]{colbourn1981testing}
C.~J. Colbourn.
\newblock On testing isomorphism of permutation graphs.
\newblock {\em Networks}, 11(1):13--21, 1981.

\bibitem[DCB13]{domahidi2013ecos}
A.~Domahidi, E.~Chu, and S.~Boyd.
\newblock {ECOS}: {A}n {SOCP} solver for embedded systems.
\newblock In {\em European Control Conference (ECC)}, pages 3071--3076. IEEE,
  2013.

\bibitem[ER63]{erdHos1963asymmetric}
R.~Erd{\H{o}}s and A.~R{\'e}nyi.
\newblock Asymmetric graphs.
\newblock {\em Acta Mathematica Hungarica}, 14(3):295--315, 1963.

\bibitem[Fru49]{frucht1949graphs}
R.~Frucht.
\newblock Graphs of degree three with a given abstract group.
\newblock {\em Canadian J. Math}, 1:365--378, 1949.

\bibitem[FS14]{fiori2014spectral}
M.~Fiori and G.~Sapiro.
\newblock On spectral properties for graph matching and graph isomorphism
  problems.
\newblock {\em arXiv preprint arXiv:1409.6806}, 2014.

\bibitem[GJ79]{garey1979computers}
M.~R. Garey and D.~S. Johnson.
\newblock Computers and intractability: {A} guide to np-completeness, 1979.

\bibitem[HW74]{hopcroft1974linear}
J.~E. Hopcroft and J.~Wong.
\newblock Linear time algorithm for isomorphism of planar graphs (preliminary
  report).
\newblock In {\em Proceedings of the sixth annual ACM symposium on Theory of
  computing}, pages 172--184. ACM, 1974.

\bibitem[Irn05]{irniger2005graph}
C.~M. Irniger.
\newblock Graph matching--filtering databases of graphs using machine learning
  techniques.
\newblock 2005.

\bibitem[JK07]{junttila2007engineering}
T.~Junttila and P.~Kaski.
\newblock Engineering an efficient canonical labeling tool for large and sparse
  graphs.
\newblock {\em Proceedings of the Ninth Workshop on Algorithm Engineering and
  Experiments (ALENEX)}, 12:135--149, 2007.

\bibitem[JK09]{Specht2013}
T.~A. Junttila and P.~Kaski.
\newblock Engineering an efficient canonical labeling tool for large and sparse
  graphs.
\newblock {http://www.tcs.hut.fi/Software/bliss/benchmarks/index.shtml}, April
  2009.

\bibitem[JK11]{junttila2011conflict}
T.~Junttila and P.~Kaski.
\newblock Conflict propagation and component recursion for canonical labeling.
\newblock In {\em Theory and Practice of Algorithms in (Computer) Systems},
  pages 151--162. Springer, 2011.

\bibitem[K{\etalchar{+}}57]{kelly1957congruence}
P.~J. Kelly et~al.
\newblock A congruence theorem for trees.
\newblock {\em Pacific J. Math}, 7(0957):961--968, 1957.

\bibitem[KBL07]{kandel2007applied}
A.~Kandel, H.~Bunke, and M.~Last.
\newblock {\em Applied graph theory in computer vision and pattern
  recognition}, volume~52.
\newblock Springer, 2007.

\bibitem[Kuh55]{kuhn1955hungarian}
H.~W. Kuhn.
\newblock The {H}ungarian method for the assignment problem.
\newblock {\em Naval research logistics quarterly}, 2(1-2):83--97, 1955.

\bibitem[LB79]{lueker1979linear}
G.~S. Lueker and K.~S. Booth.
\newblock A linear time algorithm for deciding interval graph isomorphism.
\newblock {\em Journal of the ACM (JACM)}, 26(2):183--195, 1979.

\bibitem[LFF{\etalchar{+}}14]{lyzinski2014graph}
V.~Lyzinski, D.~Fishkind, M.~Fiori, J.~T. Vogelstein, C.~E. Priebe, and
  G.~Sapiro.
\newblock Graph matching: {R}elax at your own risk.
\newblock {\em arXiv preprint arXiv:1405.3133}, 2014.

\bibitem[LPA09]{lopez2009fast}
J.~L. L{\'o}pez-Presa and A.~F. Anta.
\newblock Fast algorithm for graph isomorphism testing.
\newblock In {\em International Symposium on Experimental Algorithms}, pages
  221--232. Springer, 2009.

\bibitem[LPAC11]{lopez2011conauto}
J.~L. L{\'o}pez-Presa, A.~F. Anta, and L.~N. Chiroque.
\newblock Conauto-2.0: {F}ast isomorphism testing and automorphism group
  computation.
\newblock {\em arXiv preprint arXiv:1108.1060}, 2011.

\bibitem[Luk80]{luks1980isomorphism}
E.~M. Luks.
\newblock Isomorphism of graphs of bounded valence can be tested in polynomial
  time.
\newblock In {\em 21st Annual Symposium on Foundations of Computer Science},
  pages 42--49. IEEE, 1980.

\bibitem[SFSV03]{de2003large}
M.~De Santo, P.~Foggia, C.~Sansone, and M.~Vento.
\newblock A large database of graphs and its use for benchmarking graph
  isomorphism algorithms.
\newblock {\em Pattern Recognition Letters}, 24(8):1067--1079, 2003.

\bibitem[SU11]{scheinerman2011fractional}
E.~R. Scheinerman and D.~H. Ullman.
\newblock {\em Fractional Graph Theory: {A} Rational Approach to the Theory of
  Graphs}.
\newblock Courier Corporation, 2011.

\bibitem[Tin91]{tinhofer1991note}
G.~Tinhofer.
\newblock A note on compact graphs.
\newblock {\em Discrete Applied Mathematics}, 30(2):253--264, 1991.

\bibitem[Ume88]{umeyama1988eigendecomposition}
S.~Umeyama.
\newblock An eigendecomposition approach to weighted graph matching problems.
\newblock {\em IEEE Transactions on Pattern Analysis and Machine Intelligence},
  10(5):695--703, 1988.

\end{thebibliography}

\end{document}